\documentclass[11pt,a4paper]{article}
\pdfoutput=1
\usepackage[T1]{fontenc}
\usepackage[utf8]{inputenc}
\usepackage{lmodern}
\usepackage[english]{babel}
\usepackage[margin=2.4cm]{geometry}
\usepackage{microtype}
\usepackage{amsmath,amssymb,amsthm}
\usepackage{booktabs,array,tabularx}
\usepackage{enumitem}
\usepackage[authoryear,round]{natbib}
\usepackage{graphicx}
\usepackage{tikz}
\usetikzlibrary{calc,arrows.meta}
\usepackage{xcolor}
\usepackage[hidelinks]{hyperref}
\usepackage{xurl}

\hypersetup{
  pdftitle={One machine, three drawings: a didactic analysis of a teaching progression for the function concept from Grade 4 to Grade 10},
  pdfauthor={Vincent Tariel-Adam},
  pdfkeywords={function concept, function machine, correspondence, covariation, concreteness fading, didactics of mathematics, curriculum design}
}

\newtheorem{principle}{Design principle}
\theoremstyle{definition}
\newtheorem{hypothesis}{Hypothesis}

\colorlet{colordef}{red!60}
\definecolor{colorprop}{rgb}{0.0,0.4,0.74}
\newcommand{\imagepath}{}

\title{One machine, three drawings:\\ a didactic analysis of a teaching progression\\ for the function concept from Grade~4 to Grade~10}
\author{Vincent Tariel-Adam\\[2pt]
  \small Independent researcher; CommeUnJeu (\url{https://commeunjeu.org}), France\\
  \small \texttt{vincent.tariel@alumni.polytechnique.org}}
\date{November 2026}

\begin{document}
\maketitle

\begin{abstract}
\noindent
The function concept can be entered through a correspondence between sets, a relation that is single-valued, or through a process that turns an input into an output. This paper presents the design rationale and a critical analysis of a source progression for functions (Grades 4--10, nominal ages 9--16, with two levels on relations between quantities), from which French-curriculum and International Baccalaureate chapters published at \url{https://commeunjeu.org} are assembled. One schema is presented first and kept throughout: a function is a machine that takes an input, applies a rule and gives exactly one output. Its drawing changes with age (a robot, then a box, then the arrow $f : x \mapsto f(x)$), and each change is shown once, the two drawings side by side; after 15 the box returns only in the lessons, to explain composition and the inverse. The properties that research identifies as difficult (the uniqueness of the output, the constant function, the function without a formula, two rules for one function, preimages) are posed on the machine from Grade~4; the word \emph{relation} is kept for the link between two quantities, where covariation is taught; the set of ordered pairs comes back as the graph. The choices are related to the function machine as a cognitive root and to its authors' own warnings, to process--object accounts, to studies of concept images and of covariation, and to the French and IB curricula. A block-level audit of the dated source documents found the machine in only 2 of the 177 exercise blocks of Grades~7 to~10, no exercise asking pupils to draw a machine, and public chapters that opened on the relation; the revisions made after the audit are recorded with their counts before and after. No empirical evaluation is reported; hypotheses for one are stated.
\end{abstract}


\section{Introduction}\label{sec:intro}

The function is the concept through which school mathematics passes from computing numbers to studying dependences, and it is one of the concepts whose history weighs most on its teaching. Its definitions moved from a quantity composed of a variable and constants to an analytic expression, then to the dependence of one quantity on another, to an arbitrary law of correspondence after Dirichlet, and finally to a set of ordered pairs without any reference to a computation \citep{Kleiner1989,Sfard1991,ThompsonCarlson2017}. School definitions of the second half of the twentieth century followed the last steps: two sets and a rule that associates with each element of the first exactly one element of the second, sometimes reached from the relation as a set of pairs \citep{Markovits1986,BellJanvier1981}. Research on learners has documented the difficulties that remain: definitions recited and contradicted by the concept image \citep{VinnerDreyfus1989}, functions identified with formulas, constant and piecewise functions rejected \citep{Markovits1986,Oehrtman2008}, and prospective teachers unable to say where the requirement of a single output comes from \citep{Even1993}.

A curriculum author meets a design problem specific to this concept. The first representation chosen must carry the uniqueness of the output and the arbitrariness of the rule; it must lead to tables, formulas and graphs without being replaced by any of them; and it must not hide that two quantities can vary together. Two entries are available. The first starts from the relation, a set of pairs, and defines a function as a relation in which each first element has one partner. The second starts from a process: a machine receives an input, applies its rule, and gives one output. The second entry has been proposed as a ``cognitive root'' of the concept by the authors who also named its limits \citep{TallMcGowenDeMarois2000}.

This paper documents and analyses one answer to that problem, a progression written as bilingual (English and French) lessons and exercise sets over seven levels, with two levels on relations between quantities. Its research question is: \emph{what becomes of a concept when one schema carries it from age 9 to age 16, and where is that schema not enough?} One schema is presented first and kept: a function is a machine that takes an input, applies a rule, and gives exactly one output. The design rests on five principles:
\begin{enumerate}[label=(\roman*)]
  \item the machine is the definition, at every level; later definitions rewrite it, and its drawing lightens with age (Section~\ref{sec:p1});
  \item the properties of the concept are posed on the machine from Grade~4, with the definition as the argument (Section~\ref{sec:p2});
  \item the word \emph{relation} is kept for the link between two quantities, where covariation is taught, and the two meet in ``in terms of'' (Section~\ref{sec:p3});
  \item the set of ordered pairs is not the starting point: it comes back as the graph, the set of the points the machine gives (Section~\ref{sec:p4});
  \item $f(x)$ is read on the machine: the output that the machine $f$ gives to the input $x$ (Section~\ref{sec:p5}).
\end{enumerate}

\paragraph{Genre and scope.} The paper is a design rationale and a curriculum analysis, the companion of a study of a fraction progression by the same author \citep{TarielAdam2026Fractions}, whose method it shares. Like a hypothetical learning trajectory \citep{simon1995,clements2004}, the progression embodies learning goals, an order of activities and hypotheses about learning; the paper makes them explicit and confronts them with the literature and with the documents that implement them. Its method borrows the \emph{a priori analysis} of didactic engineering \citep{artigue1990,artigue1992} at the level of global choices (the representation, the order of meanings, the drawings) and, in Section~\ref{sec:tasks}, at the level of three tasks and their didactic variables. No classroom experiment is reported; statements about learners are taken from the cited literature or presented as design hypotheses. None of the principles is new taken alone: the function machine is classical, and its authors already recommended machines defined by a procedure rather than by a formula \citep{TallMcGowenDeMarois2000}. What the paper contributes is the articulation of one schema across seven levels and three drawings, the explicit place given to the relation and to the set of pairs, the identification of the points where the machine is not enough, and a confrontation of the stated design with the materials, which led to revisions made before this paper was written (Section~\ref{sec:revisions}). The families of functions (linear, affine, quadratic, exponential) have their own progressions and are outside the scope.

\paragraph{Grade convention.} The progression is organised in source levels labelled Grade~4 to Grade~10, nominally corresponding to ages 9 to 16, and two levels on relations labelled Grade~6 and Grade~7; the definitions of Grade~11 and of the first post-secondary year are cited to show that the definition does not change. It is a source progression, not a national curriculum: curriculum chapters are assembled from it by links (Section~\ref{sec:method}). In the French-curriculum chapters, Grade~4 serves the data chapters of \emph{CM1} and \emph{CM2} (ages 9--11) with its curve of measurements; Grades~4 to~7 and the relation levels serve \emph{cinqui\`eme} (age 12--13); the relation levels and a level on equations serve \emph{quatri\`eme}; Grades~7 to~9 serve \emph{troisi\`eme}; Grades~7 to~10 serve \emph{seconde}. In the International Baccalaureate chapters, Grade~4 serves PYP~4, Grades~4 and~5 PYP~5, the relation levels MYP~1 and MYP~2, Grades~7 to~11 MYP~5 (Table~\ref{tab:projections}). Unless stated otherwise, ``Grade~$N$'' refers to the levels of the progression.

\section{Corpus and method}\label{sec:method}

\paragraph{Corpus.} The corpus is the source documents of the function progression: a lesson and an exercise set for each of Grades 4 to 10, and for the two relation levels (Grades 6 and 7); one definition of a calculation program, written in the level on equations, is included because the progression draws it as machines in a chain. The documents are written in a mark-up in which every definition, method, example, exercise and answer is a delimited \emph{block} with a type and a stable identifier; a block of a curriculum chapter, or of a higher level, may be a \emph{link} to a block written elsewhere, which it displays without copying. Blocks are referred to by their identifier, for instance 30614. Two dated states are analysed: the state of 2~November~2026 before the revisions, on which the audit was made (commit \texttt{596ec5ecec}), and the state after the revisions of the same day (commit \texttt{e8927be52e}); Table~\ref{tab:corpus} gives their sizes.

\paragraph{Statements and counts.} As in the companion paper, three kinds of statement are made about the corpus: \emph{presence} (a block with a given content exists at a given level), \emph{order} (one block precedes another), and \emph{counts}, exhaustive over the mark-up. Counts are made by a script on the English half of each block and checked by reading; the French half is a parallel translation. A block is \emph{tagged} with a marker when its text contains it: the machine (the word, or the image of the machine), a table, a drawn graph, a formula $f(x)=\ldots$, the notation $f(\ldots)$, a calculation program, the vocabulary of relations, a correspondence presented as \emph{not} a function, the constant function, a request to \emph{draw} a machine. A tag says that a block contains the marker, not how the block uses it: every statement about a use was made by reading the block. Content counts inspect the body of the blocks a file owns; links are counted apart.

\paragraph{Projections.} The counts concern the source documents, not repeated appearances of a block in the assembled chapters. The curriculum chapters (French \emph{CM1}, \emph{CM2}, \emph{cinqui\`eme}, \emph{quatri\`eme}, \emph{troisi\`eme}, \emph{seconde}; IB PYP~4, PYP~5, MYP~1, MYP~2, MYP~5, and two chapters of MYP~5 Extended and DP~1 for composition) are described from their link tables, which list, unit by unit, the blocks they show; a second script follows each link to its source block and reports how many lesson and exercise blocks of a chapter contain the machine, and the order in which the lesson presents a relation, a machine, or a definition without a machine. The tables establish what a chapter contains and in which order, not what pupils actually experienced. The projections are counted at the state before the revisions and at the state after them (commit \texttt{d809e0071a}, which also moves the coordinate plane of the \emph{troisi\`eme} chapter just before the graphs).

\paragraph{Limits of a self-audit.} The author designed the progression and made the audit. The plan of this paper, the revision proposals and a teacher's guide drawn from the same materials were each read by a language-model reviewer, which found errors in the author's first statements (a public chapter credited with a block it does not serve; a figure whose marked point lay off its curve); the corrected statements are those reported. The source documents and the counting scripts can be consulted on request; the assembled chapters are public (Appendix~\ref{app:resources}). No public archive of the dated sources accompanies this manuscript, so the counts are not yet independently reproducible from the public chapters alone.

\begin{table}[htbp]
\centering\small
\caption{The source documents before and after the revisions of 2~November~2026 (English halves, blocks owned by each file). ``Machine'' counts blocks that contain the word or the image of the machine; ``draw'' counts blocks that match the regular expression for a request to draw a machine; it does not count every drawing task (the instruction ``Draw the chain'' of 546185 is not matched). Written exercises are answered on paper or on a slate and corrected from a photograph.}
\label{tab:corpus}
\begin{tabular}{lrrrrrr}
\toprule
 & \multicolumn{2}{c}{Lesson blocks} & \multicolumn{4}{c}{Exercise blocks} \\
\cmidrule(lr){2-3}\cmidrule(lr){4-7}
Level & all & machine & all & written & machine & draw \\
\midrule
Grade 4 & 8 & 3 & 29 & 13 & 19 & 0 \\
Grade 5 & 3 & 2 & 14 & 9 & 12 & 0 \\
Grade 6 & 3 & $1\to2$ & 13 & 13 & 6 & 0 \\
Relation, Grade 6 & 20 & $0\to1$ & 60 & 18 & 0 & 0 \\
Relation, Grade 7 & 7 & $0\to1$ & 38 & 24 & 1 & 0 \\
Grade 7 & $6\to7$ & 2 & $47\to52$ & $9\to14$ & $0\to5$ & $0\to2$ \\
Grade 8 & 7 & 0 & 19 & 14 & 1 & 0 \\
Grade 9 & 11 & $1\to2$ & $40\to41$ & $17\to18$ & $0\to1$ & 0 \\
Grade 10 & 29 & $2\to4$ & 71 & 24 & 1 & 0 \\
\bottomrule
\end{tabular}
\end{table}

\section{Theoretical background}\label{sec:background}

\subsection{From a formula to a correspondence}
The history of the concept is a succession of definitions \citep{Youschkevitch1976,Kleiner1989}. Following the account of \citet[pp.~14--16]{Sfard1991}: the word appears with Leibniz; Johann Bernoulli speaks of a quantity composed of a variable and constants, Euler first of an analytic expression and then of the dependence of one quantity on another; the successive failures lead to Dirichlet's rejection of the algorithmic approach, and to Bourbaki's purely structural definition, a set of ordered pairs without any reference to a computation. \citet{ThompsonCarlson2017}, after Boyer, date from Dirichlet a precise law of correspondence between $x$ and $y$, which may be arbitrary and very discontinuous, and note that the current definition follows Dirichlet but is written in Cartesian products and pairs. School textbooks of the second half of the twentieth century defined the function by two sets and a rule that associates with each element of the first exactly one element of the second; the British SMP course started from the relation \citep[p.~18]{Markovits1986}, and its study of functions began with relations and arrow diagrams, the definition by pairs and graphing techniques taking most of the time \citep[p.~34]{BellJanvier1981}. \citet{ThompsonCarlson2017} and \citet{Oehrtman2008} argue that the set definition solved problems of mathematicians that a beginner cannot conceive: introducing the concept by it presents the solution of a problem the learner does not have.

\subsection{Process and object}
\citet{Sfard1991} distinguishes an operational conception of a function, a process of computation or a well-defined method of passing from one system to another, from a structural conception, a set of ordered pairs. In her account, operational conceptions precede structural ones in three stages: interiorisation (computing a value of the dependent variable with a formula), condensation (thinking of the process in terms of input--output relations, which allows the learner to study functions, draw their graphs, compose two of them and find an inverse), and reification, an ontological shift visible when the learner accepts that a function need not be computable. Two side-paths are noted: the prototype (beginners imagine a linear map as soon as a function is mentioned) and a quasi-structural approach that identifies the notion with one of its representations. In the APOS tradition, a computer-based teaching that presented functions as processes, which students named and then applied by name, improved the understanding of many university students \citep{Breidenbach1992}; an \emph{action} conception, by contrast, substitutes numbers into an expression one evaluation at a time \citep{DubinskyHarel1992}, a description taken up by \citet{Oehrtman2008}: with an action view, domain and range are not thought of as sets, the inverse reduces to ``switching $x$ and $y$'' and composition to a substitution, while a process view imagines the function as an entity that receives an input, produces an output, and can be reversed. \citet{Thompson1994} adds that reaching the process conception is not trivial, and often not reached without teaching that aims at it.

\subsection{Concept image and known difficulties}
\citet{Vinner1983} built a model of learning on the distinction between the definition of a concept and the image of it that learners carry. In a study of 271 college students and 36 junior high school teachers, many definitions and even more images remained primitive, and among those who gave the modern definition, image and definition often contradicted each other \citep{VinnerDreyfus1989}. Among 152 prospective secondary teachers, many lacked a modern conception: the arbitrariness of functions escaped them, and very few could explain where the requirement of a single output comes from, or why it matters \citep{Even1993}. In a study of about 400 pupils aged 14--15, after the concept had been taught, \citet{Markovits1986} found three kinds of functions hard to accept, the constant function, the function defined piecewise and the function given by a discrete set of points; the domain and the codomain neglected; images and preimages only partly understood; a poor repertoire of examples, too close to linearity; and the passage from the graph to the formula harder than the converse. Many students, they write, ``think that a function is defined by the rule of correspondence only'' \citep[p.~22]{Markovits1986}. Among good university algebra students, according to Carlson as reported by \citet{Oehrtman2008}, few recognise a constant function as a function, many believe that a function must be given by one formula, see a piecewise function as two functions, and find it hard to accept that two different formulas define the same function.

\subsection{The function machine as a cognitive root, and its limits}
\citet{TallMcGowenDeMarois2000} propose the input--output box as a \emph{cognitive root} of the concept: a unit of knowledge that is meaningful at the beginning of learning, develops by expansion rather than reconstruction, keeps its meaning in the long term and remains useful as understanding grows. The box has the status of an object and carries the process from input to output; tables, graphs, formulas, procedures and verbal statements become ways of representing or computing the relation inside. Two boxes are ``the same'' when they give the same output for each input of the domain, whatever their inner mechanism. The paper is mainly theoretical, drawing on earlier data, and its authors name three limits of their own root. Used as a game of ``guess my rule'', the box produces ``the epistemological obstacle that all functions are given by a formula''; hence their advice of everyday machines defined by a procedure, such as the price of posting a letter according to its weight. The box has no explicit domain or codomain, and may suggest a ``natural'' domain; they advise presenting it early as an arrow from a set $A$ to a set $B$. And the graph of a real function adds the order of the reals on the axes: it is a ``function plus'', not the function alone. In the companion study of students taking non-credit algebra courses in two community colleges --- adults, not children --- two thirds of the introductory group read a machine diagram from the first day, but at the end of the semester only 41\% of the introductory group read an output on a graph and 22\% went back from the graph to the input; the authors conclude that the link between machines and graphs was weak for these students \citep{McGowenDeMaroisTall2000}. A technology-rich design for Grade~8 found that the link between tool use and conceptual development benefited from preparatory activities, from tools whose representations allow increasingly elaborate reasoning, and from alternation with paper and pencil \citep{Doorman2012}.

\subsection{Functional thinking before the function}
Children of nine and a half explored function tables by attending to the relation between inputs and outputs, and expressed it orally and in symbols \citep{Warren2006}; children of 9--10 used algebraic ideas and representations usually absent from primary school \citep{Carraher2006}; six-year-olds thought in generalised ways about functional data, against a tradition that only looks at the variation within one sequence of values \citep{Blanton2015}. In an early-algebra intervention in Grades 3--5, written responses progressed from recursive patterns to correspondence, and sometimes to covariation \citep{Stephens2017}. Students whose earlier notions had been organised around the word ``function'' generally showed a better understanding of the function as an object \citep{Ayalon2017}.

\subsection{Correspondence and covariation}
A covariational approach coordinates the passage from one value of $y$ to the next with the passage from one value of $x$ to the next, for instance down the columns of a table \citep{ConfreySmith1994,ConfreySmith1995}. \citet{Thompson1994} proposes covariation from kindergarten to the second year of university, and correspondence when its need arises, while noting that some correspondences have neither an analytic expression nor a covariational reading, such as a person's name and their social security number. \citet{Carlson2002} give a framework of the mental actions of covariational reasoning. \citet{ThompsonCarlson2017} argue that the view of a function as one input giving one output, alone, does not equip students to reason about two quantities that vary together, and recall the hypothesis that conceiving a rule of correspondence covariationally first requires a process conception of that rule. The review of \citet{Leinhardt1990} counts correspondence and covariation as the two ways of understanding a functional relation.

\subsection{Registers and graphs}
For \citet{Duval1993}, converting a piece of knowledge from one semiotic register to another is a fundamental cognitive operation, irreducible to a treatment within one register, and many learning difficulties come from a lack of coordination between registers. Graphs of situations are readily misread when their pictorial aspect conflicts with their meaning, and tables help to study how variables change \citep{BellJanvier1981}. In a study of French \emph{seconde} classes (260 pupils), \citet{CoppeDorierYavuz2007} found that tables of values and graphs may not be the most adequate entries to the notion of function, and that a conversion task can be done without understanding what a function is. The difficulty of handling functions not defined by an algebraic formula persists at the entrance to university \citep{Vandebrouck2011}; a graphic milieu in which pupils of \emph{premi\`ere} worked on functions given only by their graphs made them less dependent on classical exercises \citep{Bloch2002}.

\subsection{Concreteness fading}
The sequence in which a concrete representation is progressively replaced by more abstract ones is known as \emph{concreteness fading}, the subject of a systematic review for mathematics and science instruction \citep{Fyfe2014}. The three drawings of the progression follow that order; the paper does not claim that the review supports the particular drawings or ages chosen.

\subsection{Institutional context}
In France, the programme of \emph{troisi\`eme} of 2008 made the notion emerge on examples as a ``processus faisant correspondre \`a un nombre un autre nombre'' \citep{MEN2008Troisieme}. The resources of 2016 for the programme of \emph{cycle~4} installed the function by making pupils understand that it associates with certain numbers a unique image, relied on calculation programs, and stated that the rigorous definition of a numerical function is not an objective of secondary education \citep{Eduscol2016Fonctions}. The programme of 2020 lists vocabulary, representations (symbolic expression, table, graph, calculation program) and the notations $f(x)$ and $x \mapsto f(x)$, without a definition \citep{MEN2020Cycle4}. The programme of 2026 introduces the notion ``d'abord dans le cadre des grandeurs'', uses ``en fonction de'' from \emph{cinqui\`eme}, asks in \emph{quatri\`eme} to apply a calculation program and to find the starting number after running it backwards, studies examples of functions in \emph{troisi\`eme} ``sans \'etude g\'en\'erale de la notion de fonction'', and uses the notations $P(A)$, $p(t)$ and the arrow progressively \citep{MEN2026Cycle4}. In \emph{seconde}, the notion is consolidated ``comme exprimant la d\'ependance d'une variable par rapport \`a une autre'', the graph being the set of points whose coordinates satisfy $y=f(x)$ \citep{MEN2019Seconde,MEN2026Seconde}. In the International Baccalaureate, the Middle Years Programme guide defines a function as a relation or expression involving one or more variables, keeps this strand out of MYP~1--3, and places mappings, function notation, and linear and quadratic functions in MYP~4--5, with domain and range in the extended course and composite and inverse functions as enrichment \citep{IBO2020MYP}. The Diploma Programme guide describes functions as mappings that associate each value of the independent variable (the input) with one and only one value of the dependent variable (the output), and introduces the inverse informally as what ``reverses or undoes the effect of a function'' \citep{IBO2019DPAA}.

\section{Principle 1: the machine is the definition}\label{sec:p1}

\begin{principle}
A function is a machine that takes an input, applies a rule and gives exactly one output; the same input always gives the same output. This is the definition at every level: later definitions rewrite it with letters, sets and an arrow, never replace it by another. Its drawing lightens with age, and each change of drawing is shown once, the two drawings side by side.
\end{principle}

\subsection{One definition, rewritten}
At Grade~4 the definition is already the right one: ``A function is a rule that takes an input number and gives back exactly one output number. For each input, the function machine always gives the same output and follows the same rule every time'' (182067). A text placed before it relates the machine to the In/Out tables of number patterns, a position table being an In/Out table where the position goes in and the term comes out (3921). At Grade~7, the remark ``The Machine, Drawn as a Box'' shows the same machine twice, then gives it a name: ``Once the machine is given a name, $f$, its output is written $f(x)$, read `$f$ of $x$'\,'' (546182). At Grade~9 the definition keeps the box and adds the line: a function $f$ ``is a rule that assigns to each input value exactly one output value. It works like a machine drawn as a box [\ldots] The inputs that the machine accepts are the values of $x$ for which the rule gives an output. The same machine is also written in one line, with an arrow'' (30614). At Grade~10 the domain is ``the inputs that the machine accepts'' (3961), and at Grade~11 the rule goes from a set $X$ to a set $Y$, $f : X \to Y$ (4492). In the first post-secondary year the definition is still a rule: a map from $E$ to $F$ ``is a rule $f$ that associates with each element $x$ of $E$ one and only one element of $F$'' (22847); the graph, a part of $E \times F$, comes after, with its characterisation (22849). The progression never leaves the machine for another definition; it writes it differently.

\subsection{Three drawings}
\begin{figure}[htbp]
\centering
\begin{tikzpicture}[>=latex]
\node at (0,0) {\includegraphics[width=4cm]{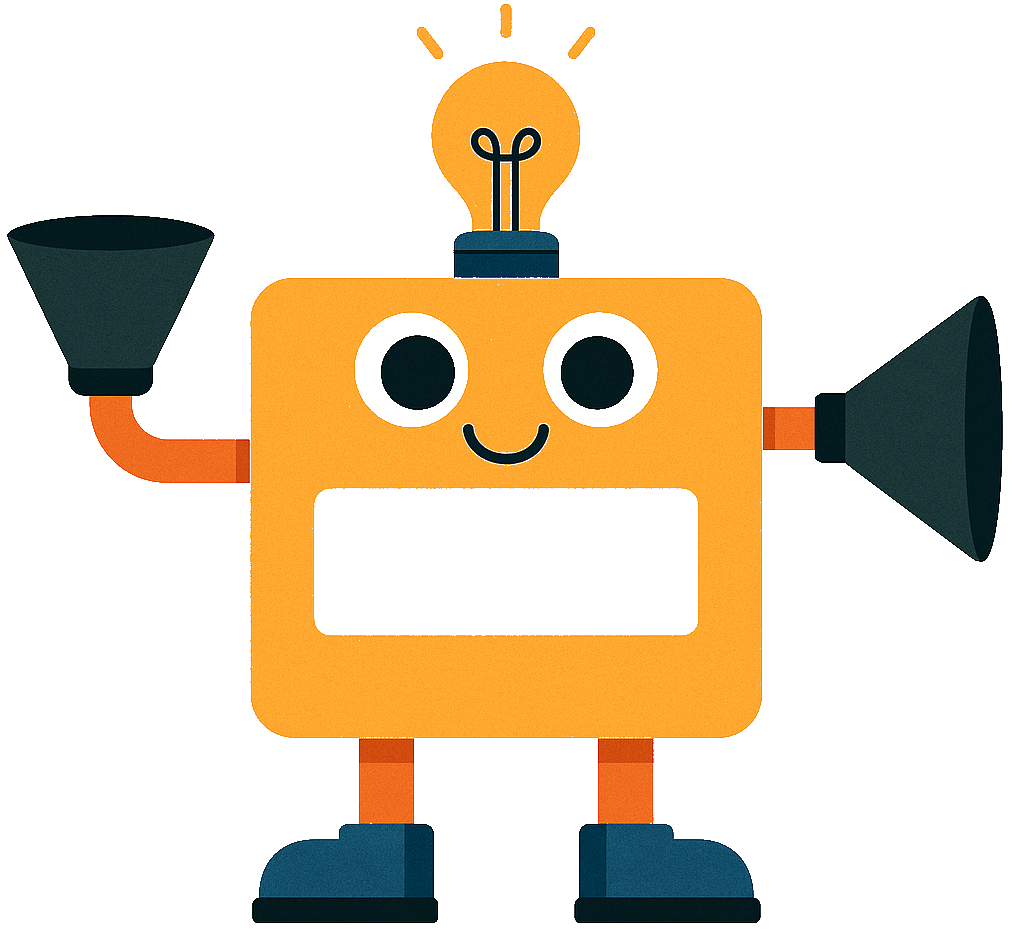}};
\node[olive] at (0,-0.4)  {rule: $\times 2$};
\draw[->,colordef,ultra thick] (-1.6,1.5)node[above]{$\textcolor{colordef}{\text{Input: } 3}$} --++(0,-0.5);
\draw[->,colorprop,ultra thick] (2,0.1) --++(0.5,0)node[right]{$\textcolor{colorprop}{\text{Output: } 6}$};
\begin{scope}[shift={(7.0,0)}]
\draw[->,colordef,ultra thick] (-1.2,0) -- (-0.1,0);
\node[colordef,left] at (-1.2,0) {$3$};
\draw[olive,thick,fill=olive!10] (0,-0.5) rectangle (2,0.5);
\node[olive] at (1,0) {$\times 2$};
\draw[->,colorprop,ultra thick] (2.1,0) -- (3.2,0);
\node[colorprop,right] at (3.2,0) {$6$};
\end{scope}
\end{tikzpicture}

\medskip
\begin{tikzpicture}[>=latex]
\draw[->,colordef,ultra thick] (-1.2,0) -- (-0.1,0);
\node[colordef,left] at (-1.2,0) {$x$};
\draw[olive,thick,fill=olive!10] (0,-0.5) rectangle (1.6,0.5);
\node[olive] at (0.8,0) {$f$};
\draw[->,colorprop,ultra thick] (1.7,0) -- (2.8,0);
\node[colorprop,right] at (2.8,0) {$f(x)$};
\node at (5,0) {is also written};
\node at (7.6,0) {$\textcolor{olive}{f} : \textcolor{colordef}{x} \longmapsto \textcolor{colorprop}{f(x)}$};
\end{tikzpicture}
\caption{The two changes of drawing, each shown once. Above, Grade~7: the robot of the earlier years and the box (block 546182). Below, Grade~9: the box and the line (block 30614). In these transition drawings the input is red, the rule olive, the output blue.}
\label{fig:drawings}
\end{figure}
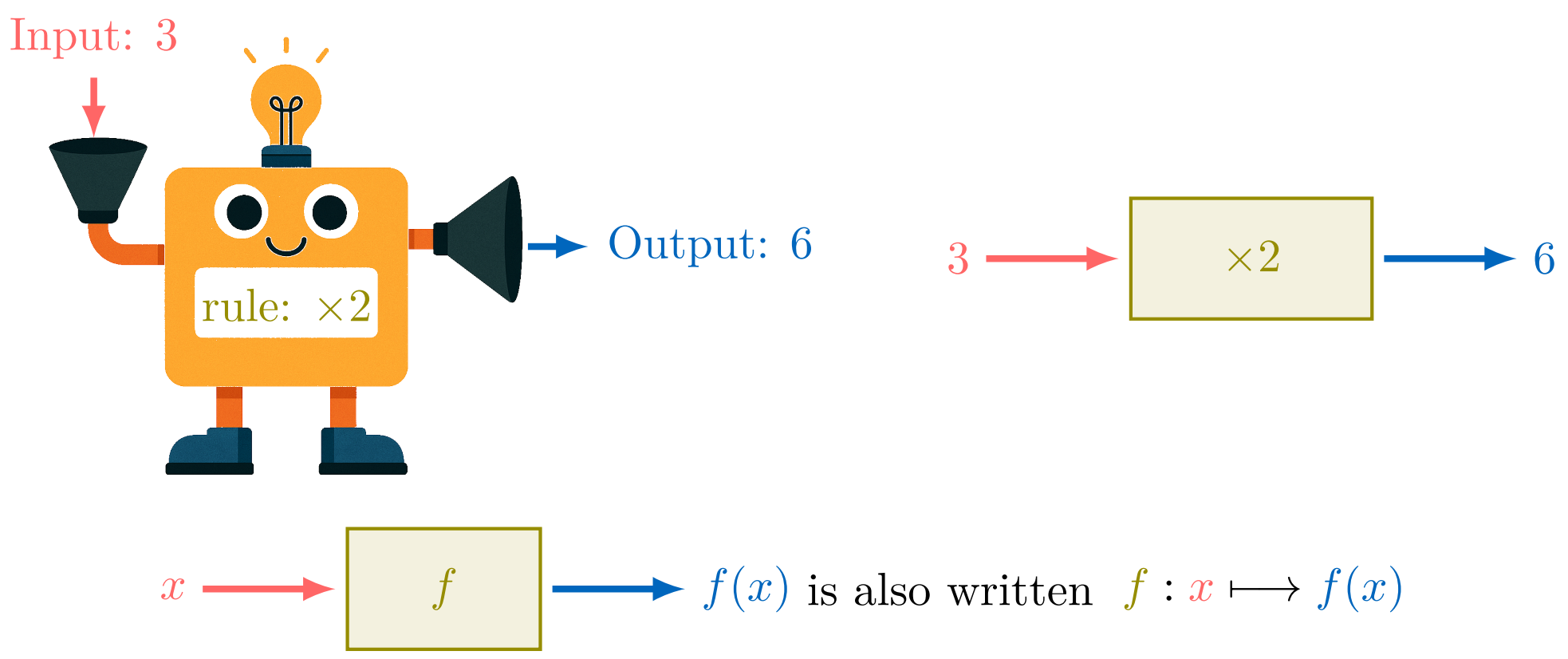

The machine of Grade~4 is a robot: the rule is written on it, the input goes in at the top, the output comes out at the side. The author judged that drawing too childish for pupils of thirteen, and the revisions of November 2026 fixed three drawings of the same machine: the robot until about 11 (Grades 4 to 6); a box, a rectangle with the rule inside and two arrows, from about 12 (Grades 7 to 9); the arrow $f : x \mapsto f(x)$ from about 15. Each change is shown once, the two drawings in one figure (Figure~\ref{fig:drawings}). After 15 the box returns only in the lessons, where it explains: two machines in a chain for composition (3969), a machine undone by the next one for the inverse (3973); no exercise asks for it.

\subsection{What the choice gains and what it postpones}
The choice gains continuity: a pupil meets at 15 the machine met at 9, which is what \citet{TallMcGowenDeMarois2000} ask of a cognitive root. It postpones the sets. The robot and the box carry no domain, the limit named by the same authors; the progression answers it in words before it answers it in symbols, ``the inputs that the machine accepts'' entering the definition at Grade~9 and naming the domain at Grade~10, the codomain waiting for Grade~11. It also postpones the codomain as part of the identity of a function: in the first post-secondary year, two maps are equal only when their domains, their codomains and their values agree (22851), and $x\mapsto x^2$ from $\mathbb{R}$ to $\mathbb{R}$ differs from $x \mapsto x^2$ from $\mathbb{R}$ to $\mathbb{R}_+$. A machine, which only shows its rule, cannot say that; Section~\ref{sec:limits} returns to it.

\section{Principle 2: the properties of the concept, posed on the machine}\label{sec:p2}

\begin{principle}
The properties that make the concept difficult --- one input never has two outputs, two inputs may share an output, a machine may answer the same number to every input, a rule need not be a formula, two rules may define the same function, a machine may be run backwards --- are posed on the machine from Grade~4, in written exercises whose argument is the definition.
\end{principle}

The written exercises of Grades 4 and 5, answered on paper or on a slate and corrected from a photograph, already carry the difficulties that the literature describes, each with the definition as the argument. A machine that answered $10$ yesterday and $12$ today for the same input $5$ cannot be a function machine (291812). A machine that answers $7$ to every input is one, against a pupil's claim that ``it doesn't compute anything'' (291832): the constant function of \citet{Markovits1986} is met at nine or ten, as a machine. A mystery machine asks for a rule that fits \emph{every} column of a table, because ``one result can suggest a rule; the results can reject a rule, never prove it'' (291818, correction 291819; 291826); the corrections of these written exercises accept any output consistent with the rule the pupil found, which avoids treating a finite set of results as determining a unique rule, but does not by itself counter the belief that every function has a formula. A machine run backwards finds the input of a known output, then checks it forwards (291830). And two machines with different rules, ``multiply by $2$'' and ``add the input to itself'', are the same machine because they give the same output for every input (291836): the criterion of sameness of \citet{TallMcGowenDeMarois2000}, stated to pupils of ten.

At Grade~7 a written group, added by the revisions, asks the pupil to draw the machine (546183--546191): $f(x) = 3x+1$ drawn as two machines in a chain, ``multiply by $3$'' then ``add $1$''; two machines exchanged, which change the function ($5(x-2)$ against $5x-2$); a chain run backwards from the output $22$, the last machine undone first; three functions $u(x) = 2(x+3)$, $v(x) = 2x+6$, $w(x) = 2x+3$, where tests on three inputs only suggest that $u$ and $v$ are equal and the distributive law proves it for every input; and a post office machine whose rule is a table of steps, the price jumping at $20$~g and at $100$~g, which is a function although no single formula $a \times m + b$ writes it --- the everyday machine that \citet{TallMcGowenDeMarois2000} recommend, and a piecewise function in the sense of \citet{Markovits1986}. At Grades~7 and~8 other written exercises ask for reasons on arrow diagrams between pupils and colours, where the whole test is ``one arrow out of every input'' (291868, correction 291869), on drawings tested with vertical lines (291882), on a pupil's claim that a function may not give the same output twice, refuted by $(-2)^2 = 2^2$ (291886), on two tables of which one is constant (291888), and on situations of the class, one of which fails because the whole numbers from $13$ to $31$ have no month (291898).

\section{Principle 3: the relation, for two quantities that vary together}\label{sec:p3}

\begin{principle}
The word \emph{relation} is kept for the link between two quantities: a table of measurements, a graph of a journey, a conversion, a formula. When each value of the first quantity determines exactly one value of the second, the second is expressed \emph{in terms of} the first, and the dependence is drawn as the machine; the relation adds what the machine alone does not show, the covariation.
\end{principle}

\begin{figure}[htbp]
\centering
\begin{tikzpicture}[>=latex]
\draw[->,colordef,ultra thick] (-1.2,0) -- (-0.1,0);
\node[colordef,left] at (-1.2,0) {side $c$};
\draw[olive,thick,fill=olive!10] (0,-0.5) rectangle (2,0.5);
\node[olive] at (1,0) {$\times 4$};
\draw[->,colorprop,ultra thick] (2.1,0) -- (3.2,0);
\node[colorprop,right] at (3.2,0) {perimeter $P$};
\end{tikzpicture}
\caption{``In terms of'', Relation Grade~7 (block 331804): the relation between the side and the perimeter of a square, drawn as the machine, before the sentence of covariation.}
\label{fig:interms}
\end{figure}
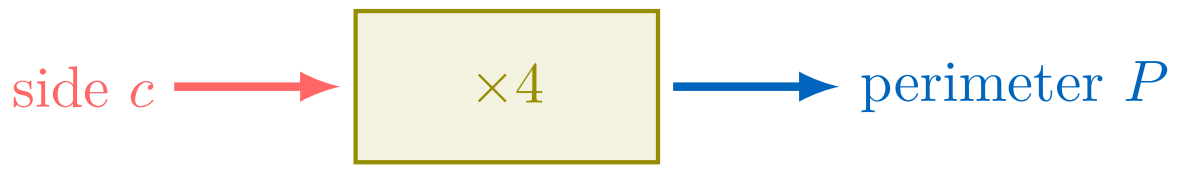

The two relation levels work on situations. At Grade~6: ``A relation describes how two variables (quantities) are connected'' (201930); tables and line graphs of measurements; distance-time graphs, where the time goes across because the distance depends on it, and each point is read as a sentence (444316); conversion graphs, where a conversion works like a machine, euros in and dollars out under the rule ``multiply by $1.2$'', the graph read backwards running the machine backwards (444323). At Grade~7 the two meet. ``Linked'' is not ``determined'': the shoe sizes and the heights of five pupils are connected, but the size $37$ appears twice with two heights, so no rule computes the height from the size (292062). When each value of the first quantity determines exactly one value of the second, the second is expressed in terms of the first, and the definition draws it (Figure~\ref{fig:interms}), then adds the covariation: ``when the side increases by $1$ cm, the perimeter increases by $4$ cm, whatever the side'' (331804). The method that follows writes the same dependence as a formula, $P = 4 \times c$, in functional notation, $P(3) = 12$, in a table and in a graph (331805).

The relation is also written, at Grade~6, as a set of ordered pairs, with a counterexample that is not a function, ``is a factor of'' from $\{2, 3\}$ to $\{4, 6, 9\}$, where $2$ has two partners (450209). That remark is the only place, before the graph, where the progression writes a relation as a set of pairs, and it is placed after the machine wherever it is served (Section~\ref{sec:progression}). The design hypothesis is that covariation, which \citet{ThompsonCarlson2017} find missing from a view of functions as ``one input, one output'', is taught where it is natural, between quantities, and joins the machine at ``in terms of''. Whether the joint holds is not shown by the documents: covariation is said once in the definition of ``in terms of'' and practised in the relation levels, not in the function levels (Section~\ref{sec:discussion}).

\section{Principle 4: the set of pairs comes back as the graph}\label{sec:p4}

\begin{principle}
The set of ordered pairs is not the starting point. It comes back as the graph: the points whose coordinates are an input and the output the machine gives to it. The machine gives the points, not the line drawn between them.
\end{principle}

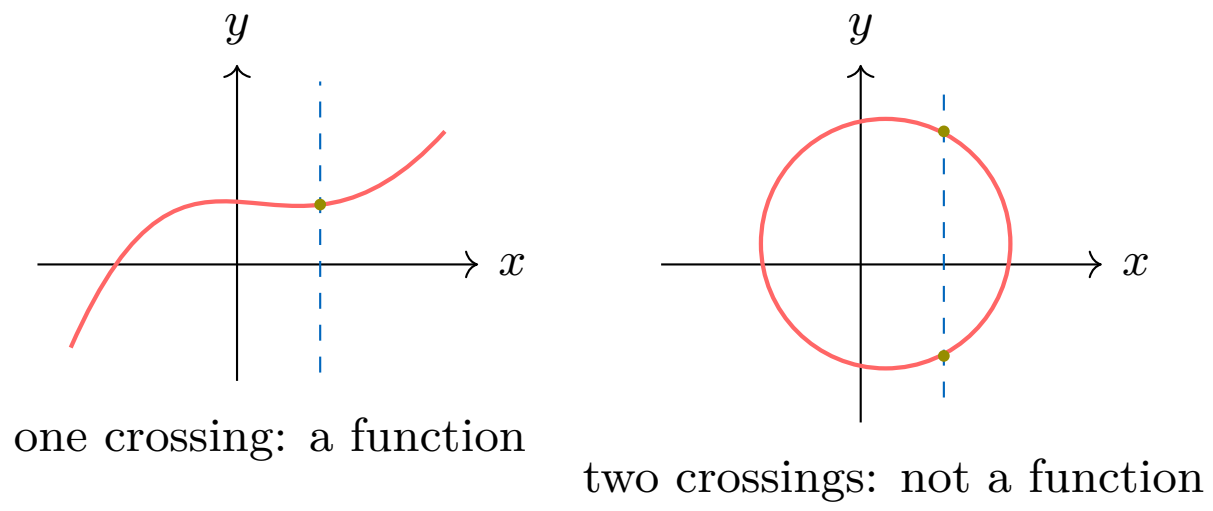
\begin{figure}[htbp]
\centering
\begin{tikzpicture}[scale=0.55]
\begin{scope}
\draw[->](-2.4,0)--(2.9,0)node[right]{\footnotesize $x$};\draw[->](0,-1.4)--(0,2.4)node[above]{\footnotesize $y$};
\draw[colordef,thick](-2,-1)..controls(-0.5,2.4)and(0.5,-0.6)..(2.5,1.6);
\draw[dashed,colorprop](1,-1.3)--(1,2.2);
\fill[olive](1,0.72)circle(2pt);
\node[font=\footnotesize]at(0.4,-2.1){one crossing: a function};
\end{scope}
\begin{scope}[xshift=7.5cm]
\draw[->](-2.4,0)--(2.9,0)node[right]{\footnotesize $x$};\draw[->](0,-1.9)--(0,2.4)node[above]{\footnotesize $y$};
\draw[colordef,thick](0.3,0.25)circle(1.5);
\draw[dashed,colorprop](1,-1.6)--(1,2.2);
\fill[olive](1,1.6)circle(2pt);\fill[olive](1,-1.1)circle(2pt);
\node[font=\footnotesize]at(0.4,-2.6){two crossings: not a function};
\end{scope}
\end{tikzpicture}
\caption{The vertical line test, Grade~8 (block 291608): the definition read on the graph.}
\label{fig:vertical}
\end{figure}

The first curve is met at Grade~4, before any machine, in the data chapters: each measurement becomes one point that carries two readings, the points are joined in order, and the lesson says what the curve adds --- the shape of the change --- and what it does not: ``the curve joins the measurements; it does not invent the ones that were never taken'' (357157). A written exercise sets the trap: a pupil reads $20\,^\circ$C on the line halfway between August and September and thinks it was measured (357169). The sentence returns with the machine and with formulas. At Grade~6, ``Only the points are results of the machine; the segments are drawn to help the eye follow them from one point to the next'' (202177); at Grade~7, between two computed points, the shape of the line ``remains a conjecture; computing one more value with the formula can rule it out, never prove it'' (422895); at Grade~8 the graph of a function is the set of all points $(x, f(x))$, and ``a table of values gives only some of these points: joining them is a conjecture until the rule gives the values in between'' (30470). The vertical line test reads the definition on the graph: the points of a vertical line all share the same input, so two crossings would give one input two outputs (291608, Figure~\ref{fig:vertical}). In the first post-secondary year, the graph becomes a part of $E\times F$ characterised by $\forall x\in E,\ \exists!\, y\in F,\ (x,y)\in\Gamma$ (22849): the set of pairs arrives last, as a property of the graph of a map.

For \citet{TallMcGowenDeMarois2000}, the graph of a real function is a ``function plus'': it adds the order of the reals on the axes. The progression does not say it in those words, but the repeated sentence on joining points is where the addition shows: the line between two points is a claim about inputs that were not computed.

\section{Principle 5: $f(x)$ read on the machine}\label{sec:p5}

\begin{principle}
The notation $f(x)$ is a reading of the machine: $f$ is the name of the machine, $x$ its input, $f(x)$ the output it gives to $x$. It is introduced on the machine, never as a formula that would come first.
\end{principle}

The revised blocks present the notation on the box. At Grade~7 the remark that draws the box names it $f$ and writes its output $f(x) = 2x$ (546182); the example that follows turns a calculation program, ``Choose a number. Multiply the chosen number by $3$. Add $1$ to the result'', into a function: the program is run on $x$ step by step, $x \to 3x \to 3x+1$, so $f(x) = 3x + 1$, and checked on a number (331826). At Grade~9 the opening problem goes from the box ``multiply by $2$'' to its table, then to $f(\text{input}) = \text{output}$, the brackets indicating ``the action of the function $f$ on the input'', then to $f(x) = 2x$ (30613); the definition names the parts: $f$ the function (the rule), $x$ the input variable, $f(x)$ the output, ``the image of $x$ under $f$'', and $x$ an antecedent, also called a preimage, of $f(x)$ (30614). In the relation levels, where letters of quantities become names of functions, the method warns that $P(3)$ does not mean ``$P$ times $3$'' but ``the perimeter when the side is $3$'' (331805). The design hypothesis is that a pupil who reads $f(x)$ as the output of $f$ for the input $x$ computes $f(x + a)$ by passing $x + a$ through the machine, not by adding $a$ at the end of the expression --- the error that, according to Carlson as reported by \citet{Oehrtman2008}, 43\% of good university algebra students made.

\section{The progression, level by level}\label{sec:progression}

\paragraph{Grades 4 and 5: the robot, the table, the first curve.} The machine, its definition and its In/Out tables (Grade~4), then the table of values as the machine's notebook (Grade~5, 193109), and the written exercises of Section~\ref{sec:p2}. In the same years, a data chapter places one point for each measurement and says what the curve does not show (357155--357157). In the French curriculum, which has no functions in \emph{cycle~3}, only the curve is served, in \emph{CM1} and \emph{CM2}; the IB serves the machine in PYP~4 and PYP~5, and the curve in PYP~4.

\paragraph{Grade 6 and the relation levels: quantities that vary together.} The machine's pairs placed as points (Grade~6); relations, tables, line graphs, journeys and conversions (Relation, Grade~6); ``in terms of'', the representations of a dependence, when to join points (Relation, Grade~7). The relation levels serve MYP~1 and MYP~2, and the French \emph{cinqui\`eme} and \emph{quatri\`eme}, which also serve the calculation program of the level on equations: ``each arrow works as a small machine with its own rule, and the whole program is these machines in a chain'', run backwards ``with the inverse operations, in the reverse order'' (299656).

\paragraph{Grade 7: the box and its name.} The passage from the robot to the box, the name $f$ and $f(x)$, the calculation program turned into a function, and the written group that asks pupils to draw functions as machines (Section~\ref{sec:p2}). It serves the French \emph{cinqui\`eme}, where the chapter now opens on the machine, and, through the later levels, \emph{troisi\`eme}, \emph{seconde} and MYP~5.

\paragraph{Grades 8 and 9: the graph, the image, the preimages.} The graph as the set of points $(x, f(x))$, the vertical line test, reading an image on a graph (Grade~8); the box and the line side by side, images and preimages read on a graph, in a table and computed, the square function and its preimages (Grade~9). They serve \emph{troisi\`eme}, \emph{seconde} and MYP~5.

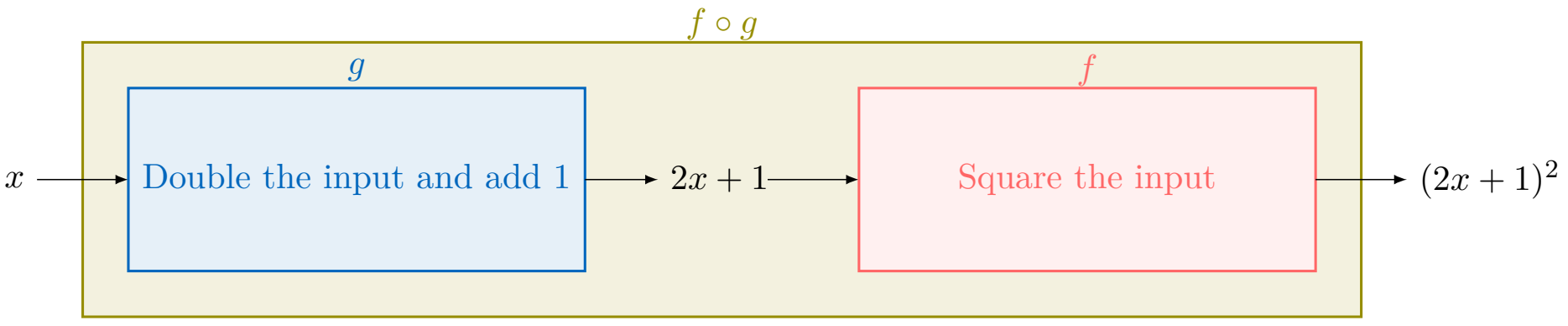
\begin{figure}[htbp]
\centering
\scalebox{0.72}{%
\begin{tikzpicture}[>=latex]
\draw[olive,thick,fill=olive!10] (-0.5,-1.5) rectangle ++(14,3);
\node[olive] at (6.5,1.7){$f \circ g$};
\draw[->] (-1,0)node[left] {$x$} -- (0,0) ;
\draw[colorprop,thick,fill=colorprop!10] (0,-1) rectangle (5,1) node[pos=.5] {Double the input and add 1};
\node[colorprop] at (2.5,1.2){$g$};
\draw[->] (5,0) -- (5.8,0)node[right] {$2x+1$} ;
\begin{scope}[shift={(8,0)}]
\draw[->] (-1,0) -- (0,0) ;
\draw[colordef,thick,fill=colordef!10] (0,-1) rectangle (5,1) node[pos=.5] {Square the input};
\node[colordef] at (2.5,1.2){$f$};
\draw[->] (5,0) -- (6,0)node[right] {$(2x+1)^2$} ;
\end{scope}
\end{tikzpicture}}
\caption{Composition, Grade~10 (block 3969): after 15 the box returns only in the lessons, where it explains. Scaled from the lesson figure.}
\label{fig:composition}
\end{figure}

\paragraph{Grade 10: the inputs the machine accepts, and the box again.} The domain and the range, even and odd functions, operations on functions, composition drawn as two machines in a chain (Figure~\ref{fig:composition}), and the inverse introduced by a machine undone by the next one (3973); for an inverse to exist, ``each output $y$ in the range of $f$ must correspond to exactly one input $x$'' (3974). The domain serves \emph{seconde} and MYP~5; composition serves MYP~5 Extended, DP~1 and the French \emph{terminale}, the inverse MYP~5 Extended and DP~1.

Table~\ref{tab:projections} gives, for each public chapter, the number of lesson and exercise blocks that contain the machine, before and after the revisions, and the order in which its lesson presents a relation, a machine and a definition without a machine.

\begin{table}[htbp]
\centering\small
\caption{The public chapters before and after the revisions: lesson and exercise blocks that contain the machine (over all linked blocks, answers and videos excluded), and the order of the lesson's tagged blocks (M: a block that contains the machine; R: a definition of a relation; F: a definition of a function without the machine). Sources are levels of the progression (R: relation levels; Eq.: the level on equations). Before: commit \texttt{596ec5ecec}; after: \texttt{d809e0071a}.}
\label{tab:projections}
\footnotesize\setlength{\tabcolsep}{4pt}
\begin{tabularx}{\textwidth}{l>{\raggedright\arraybackslash}Xllll}
\toprule
Year & Chapter & Sources & Lesson & Exercises & Order \\
\midrule
CM1 & Organisation et gestion de donn\'ees & 4 & 0/19 & 0/61 & -- \\
CM2 & Organisation et gestion de donn\'ees & 4 & 0/15 & 0/42 & -- \\
5e & Fonctions & 4--7, R6, R7 & $3\to4$/18 & $9/59\to12/62$ & RMMM $\to$ MMMRM \\
4e & Fonctions & R6, R7, Eq.~7 & $0\to2$/26 & 0/60 & R $\to$ RMM \\
3e & Fonctions & 7, 8, 9 & $2\to3$/42 & $0/98\to6/109$ & MFM $\to$ MMM \\
2de & G\'en\'eralit\'es sur les fonctions & 7--10 & $0\to2$/23 & 0/60 & F $\to$ MM \\
\midrule
PYP 4 & Patterns and functions & 4 & 3/10 & 19/52 & MMM \\
PYP 5 & Number sequences and rules & 4, 5 & 4/12 & 12/51 & MMMM \\
MYP 1 & Line Graphs & R6, 7 & $0\to1$/15 & 0/59 & -- $\to$ M \\
MYP 2 & Coordinates and Straight-Line Graphs & R6, R7 & $0\to1$/10 & 0/48 & -- $\to$ M \\
MYP 5 & Functions & 7--11, R6 & $1\to3$/39 & $0/97\to1/98$ & RMF $\to$ MMRM \\
MYP 5 Ext. & Further Functions & 10, 11 & $2\to3$/27 & 1/84 & MM $\to$ MMM \\
DP 1 & Functions (AA SL) & 10, 11 & 2/43 & $0/68\to0/72$ & MM $\to$ FMM \\
\bottomrule
\end{tabularx}
\par\smallskip\raggedright\footnotesize In DP~1 the four new exercises come from separate work on exam-style questions, and the definition of Grade~11 (4492) lost its robot, hence F. No MYP~3 or MYP~4 chapter serves these documents.
\end{table}

\section{Three tasks, and where the machine is not enough}\label{sec:tasks}

\subsection{Two rules, one function}\label{sec:task-equal}
The written exercise asks: three functions are defined by $u(x) = 2(x+3)$, $v(x) = 2x+6$ and $w(x) = 2x+3$; draw each as machines in a chain; compute them for the inputs $0$, $1$ and $5$; identify the two that seem to be the same, then prove it for every input; give an input that shows that $w$ differs (546189, Grade~7, served in \emph{cinqui\`eme} and \emph{troisi\`eme}). Its didactic variables are the drawing, which suggests a difference --- $u$ is the chain $+3$ then $\times 2$, $v$ the chain $\times 2$ then $+6$ --- the number of inputs tested, which can only suggest the equality, and the third function, which differs from both and is refuted by one input. The criteria of correction penalise the conclusion drawn from the three tests alone, and $u$ and $v$ declared different because their chains differ. The task puts in front of the pupil the criterion of \citet{TallMcGowenDeMarois2000}, the same output for each input whatever the mechanism, and the difficulty reported by \citet{Oehrtman2008}, two different formulas for one function. It also shows the limit of that criterion at school: ``each input'' means every real number, a domain the machine does not draw and the exercise does not name.

\subsection{Preimages, and why running backwards is not always a function}\label{sec:task-inverse}
\begin{figure}[htbp]
\centering
\begin{tikzpicture}[xscale=0.8,yscale=0.4]
\draw[-stealth] (-3.6,0)--(3.6,0) node[right]{\small $x$};
\draw[-stealth] (0,-0.6)--(0,9.9) node[above]{\small $y$};
\foreach \x in {-3,-2,-1,1,2,3}{\draw (\x,0.15)--(\x,-0.15) node[below]{\tiny $\x$};}
\foreach \y in {1,4,9}{\draw (0.08,\y)--(-0.08,\y) node[above left]{\small $\y$};}
\draw[colorprop,thick] plot[smooth] coordinates {(-3.1,9.61) (-3,9) (-2.5,6.25) (-2,4) (-1.5,2.25) (-1,1) (-0.5,0.25) (0,0) (0.5,0.25) (1,1) (1.5,2.25) (2,4) (2.5,6.25) (3,9) (3.1,9.61)};
\draw[olive,thick] (-3.6,9)--(3.6,9) node[right]{\small $y = 9$};
\draw[colordef,dashed] (-3,9)--(-3,0);
\draw[colordef,dashed] (3,9)--(3,0);
\fill[colordef] (-3,9) circle (2.4pt);
\fill[colordef] (3,9) circle (2.4pt);
\end{tikzpicture}
\caption{The square machine run backwards, Grade~9 (correction of exercise 546193): the line $y=9$ meets the curve above the two inputs $-3$ and $3$.}
\label{fig:square}
\end{figure}
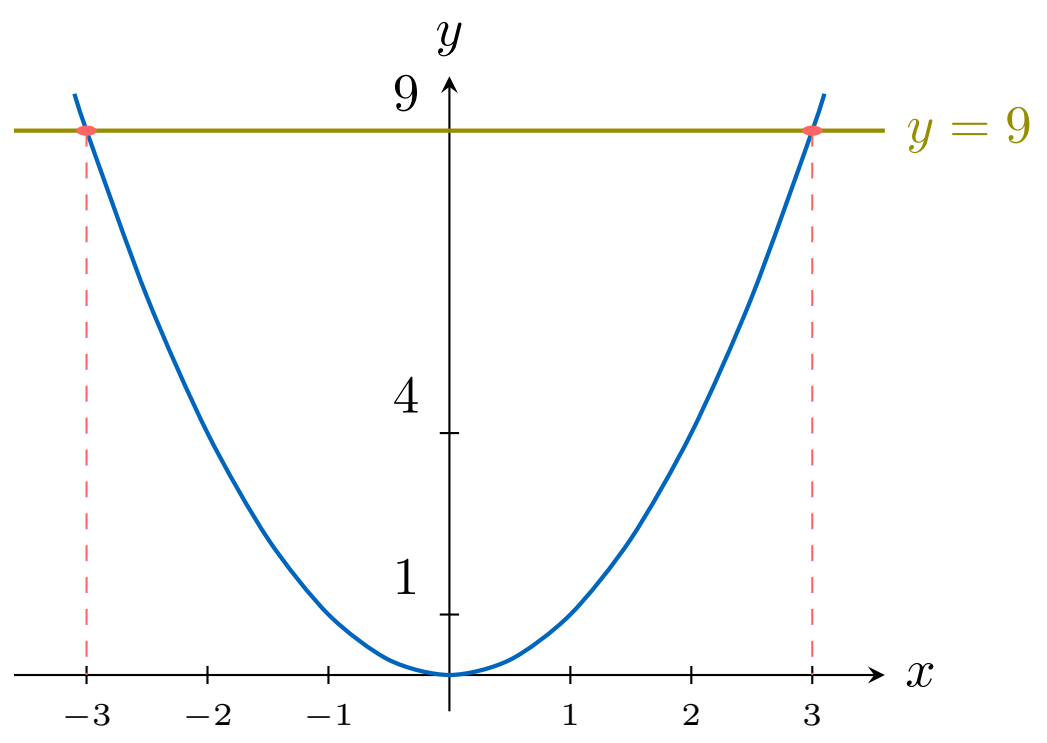

The written exercise asks: the square machine has the rule ``square the input'' and accepts every number; pass $4$ and $-4$ through it; run it backwards to find the inputs that give $9$, then $0$, then $-1$; draw the curve for $-3\le x\le 3$ and the line $y = 9$, and read the inputs that give $9$; decide whether the square machine can be run backwards like the machine ``add $3$'', with exactly one input for each output (546193, Grade~9, served in \emph{troisi\`eme} and MYP~5). Its variables are the output, positive, zero or negative, which gives two, one or no input --- the three cases of the lesson on the square function (422688) --- and the register: the curve and the horizontal line show the two inputs that knowing only that the square root of $9$ is $3$ hides (Figure~\ref{fig:square}). The last question separates finding preimages from inverting: backwards, $9$ would have two outputs, so running the machine backwards does not give a function; the square-root machine undoes the square machine only for inputs positive or zero. The criteria penalise only $3$ for the output $9$, and ``yes, the square root'' without the restriction. A further answer anticipated in this design is \emph{not a function, since two inputs give $9$}, which confuses the requirement of one output with that of one input; a pupil's claim to that effect is tested in the same chapter (291886). Reversing a reasoning was the hardest step for the adults of \citet{McGowenDeMaroisTall2000}, on a machine and still more on a graph.

\subsection{From the machine to the graph}\label{sec:task-graph}
The written exercise of Grade~4 asks whether a pupil who reads $20\,^\circ$C on the line halfway between the August point and the September point may say that it was $20\,^\circ$C between August and September; the expected answer says that one temperature was recorded per month, that only the points are measurements, and that the line only guides the eye (357169, served in \emph{CM1} and PYP~4). Its variable is the nature of the starting quantity: a month taken as a time that runs, or as a label of one measurement. The same question returns at Grade~7 in its general form --- join the points only if the starting quantity can take every value in between, and even then the shape remains a conjecture (422895) --- and at Grade~8 for functions given by a rule (30470). The task addresses two findings: graphs of situations are misread when their pictorial aspect conflicts with their meaning \citep{BellJanvier1981}, and a conversion from a graph to a table can be performed without understanding what a function is \citep{CoppeDorierYavuz2007}. It does not, by itself, build the link between machines and graphs that \citet{McGowenDeMaroisTall2000} found weak; that link is carried by the definition of the graph as the set of the points the machine gives, and by the vertical line test.

\subsection{Where the machine is not enough}\label{sec:limits}
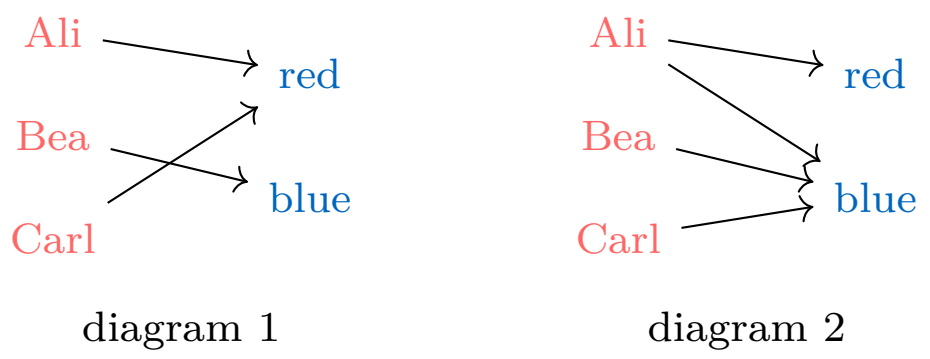
\begin{figure}[htbp]
\centering
\begin{tikzpicture}[scale=0.68]
\begin{scope}
\node[colordef] (a) at (0,2) {\scriptsize Ali};\node[colordef] (b) at (0,1) {\scriptsize Bea};\node[colordef] (c) at (0,0) {\scriptsize Carl};
\node[colorprop] (r) at (2.5,1.6) {\scriptsize red};\node[colorprop] (v) at (2.5,0.4) {\scriptsize blue};
\draw[->] (a)--(r);\draw[->] (b)--(v);\draw[->] (c)--(r);
\node[font=\scriptsize] at (1.25,-0.9) {diagram 1};
\end{scope}
\begin{scope}[xshift=5.5cm]
\node[colordef] (a) at (0,2) {\scriptsize Ali};\node[colordef] (b) at (0,1) {\scriptsize Bea};\node[colordef] (c) at (0,0) {\scriptsize Carl};
\node[colorprop] (r) at (2.5,1.6) {\scriptsize red};\node[colorprop] (v) at (2.5,0.4) {\scriptsize blue};
\draw[->] (a)--(r);\draw[->] (a)--(v);\draw[->] (b)--(v);\draw[->] (c)--(v);
\node[font=\scriptsize] at (1.25,-0.9) {diagram 2};
\end{scope}
\end{tikzpicture}
\caption{A function that is not numerical, Grade~7 (exercise 291868): diagram~1 is a function, two pupils sharing red included; diagram~2 is not, Ali sends two arrows.}
\label{fig:arrows}
\end{figure}

\paragraph{Functions that are not numerical.} The machine of Grade~4 ``takes an input number'' and gives an output number. The functions between other sets appear only in the written exercises of Grades~7 and~8 --- pupils sent to colours on arrow diagrams (Figure~\ref{fig:arrows}), the ages of the pupils of a class (291868, 291898) --- and in the definition of Grade~11, from a set $X$ to a set $Y$. \citet{Thompson1994}'s example of a correspondence with neither formula nor covariation, a name and an identification number, has no counterpart in the lessons.

\paragraph{Functions given only by a graph.} A curve of measurements has no rule to put in a box. The progression treats such curves as given points and says what was not measured; it does not ask pupils to build functions under constraints from their graphs, the work that \citet{Bloch2002} describes and that \citet{Markovits1986} recommend (a function through given points, a function increasing then constant).

\paragraph{The domain and the codomain.} The box has neither, as its authors warned \citep{TallMcGowenDeMarois2000}. The progression gives the domain in words from Grade~9 (``the inputs that the machine accepts''), asks for it in the post office task and in the situation where the numbers $13$ to $31$ have no month, and names it at Grade~10; the codomain, and its part in the equality of two functions (22851), wait for Grade~11 and the post-secondary year. The advice of presenting the machine early as an arrow from a set $A$ to a set $B$ is not followed: the arrow arrives at 15, the sets after.

\paragraph{The action view.} A machine used only to compute one output at a time can fix the action conception \citep{Thompson1994,Oehrtman2008}. The tasks of chains, sameness and running backwards aim at the process view; but from Grade~7 many exercise blocks contain symbolic expressions (in the English halves of Grade~7, 36 of the 52 exercise blocks contain a formula $f(x)=\ldots$ after the revisions), and the drawn machine is asked for in the drawing group only. Two further limits follow from the literature. \citet{Sfard1991} sees reification in the recognition that a function need not be computable; a machine with a rule inside suggests the opposite, and the progression's only answers are the table of steps and the arrow diagrams. And the machine carries a pointwise perspective, one input and its output, while the study of functions also needs global and local ones --- variation over an interval, parity, limits --- which \citet{Vandebrouck2011} asks teachers to make explicit; parity and operations on functions are taught at Grade~10 without the machine.

\section{Anticipated difficulties}\label{sec:misconceptions}

Table~\ref{tab:difficulties} lists documented difficulties and anticipated errors, and the blocks of the progression that answer them. The answers are design choices, not tested remedies.

\begin{table}[htbp]
\centering\small
\caption{Documented difficulties and anticipated errors, and the blocks that answer them.}
\label{tab:difficulties}
\begin{tabularx}{\textwidth}{>{\raggedright\arraybackslash}p{4.1cm}>{\raggedright\arraybackslash}p{3.3cm}X}
\toprule
Difficulty & Literature & Answer in the progression (level) \\
\midrule
A function is a formula & \citet{TallMcGowenDeMarois2000,Markovits1986,Oehrtman2008} & ``a rule that fits'' in mystery machines (4); the post office machine, a table of steps (7) \\
A constant function is not a function & \citet{Markovits1986,Oehrtman2008} & the machine that answers $7$ to every input (5); the constant table (8) \\
The single-output condition misapplied & \citet[p.~20]{Markovits1986} & two pupils sharing a colour (7); $(-2)^2=2^2$ (8); two letters with one price (7); the square machine backwards (9) \\
Two formulas, two functions & \citet{Oehrtman2008} & ``multiply by $2$'' and ``add the input to itself'' (5); $2(x+3)$ and $2x+6$ (7) \\
The line between points was measured & design hypothesis (related: \citealt{BellJanvier1981}) & the curve of measurements and Sami's reading (4); the conjecture between computed points (R7, 8) \\
A graph of a journey is a picture of the route & \citet{BellJanvier1981} & each point read as a sentence; distance from the start is not distance travelled (R6) \\
Images and preimages confused & \citet{Markovits1986} & the words placed on the machine: the image is the output, a preimage an input (9) \\
$f(x+a)$ computed by adding $a$ & \citet{Oehrtman2008} & $f(x)$ read on the box (7, 9); $P(3)$ is not ``$P$ times $3$'' (R7) \\
The domain neglected & \citet{Markovits1986,TallMcGowenDeMarois2000} & ``the inputs that the machine accepts'' (9, 10); numbers with no month (8) \\
A program undone in the same order & design hypothesis (the 2026 programme asks to run it backwards) & the inverse operations in the reverse order (Eq.~7); the chain from $22$ back to $7$ (7) \\
\bottomrule
\end{tabularx}
\end{table}

\section{Discussion}\label{sec:discussion}

\subsection{Position with respect to existing recommendations}
The progression follows the cognitive root of \citet{TallMcGowenDeMarois2000} and answers two of the limits its authors name: it gives machines defined by a procedure, a table of steps included, and states the domain in words before symbols. It does not follow their recommendation to introduce an arrow between explicit sets early; their warning about the additional structure of graphs is discussed in Section~\ref{sec:p4}. It diverges from \citet{Thompson1994} and \citet{ThompsonCarlson2017}, who place covariation first: the progression puts correspondence first, as the machine, and covariation in the relation levels, joined at ``in terms of''. It is consistent with the French texts on two points --- the unique image and the calculation program \citep{Eduscol2016Fonctions}, the program applied and run backwards in \emph{quatri\`eme} \citep{MEN2026Cycle4} --- and departs from the order of the 2026 programme, which introduces the notion ``d'abord dans le cadre des grandeurs'': the \emph{cinqui\`eme} chapter now opens on the machine, although on a dependence between quantities (the number of wheels is the number of bicycles multiplied by $2$, 3921), and its relation sections follow. In the IB, the MYP guide keeps the function strand out of MYP~1--3 \citep{IBO2020MYP}; the progression's MYP~1 and MYP~2 chapters carry relations, and the box appears there only for conversions and ``in terms of''.

\subsection{What the audit revealed}
On the state of 2~November~2026 before the revisions, the audit found six gaps between the stated design and the documents. (i) The machine almost disappeared from the exercise blocks once $f(x)$ was introduced: it was contained in 19 of the 29 exercise blocks of Grade~4 and 12 of the 14 of Grade~5, but in 0 of 47 at Grade~7, 1 of 19 at Grade~8, 0 of 40 at Grade~9 and 1 of 71 at Grade~10, the occurrences at Grades 8 and 10 being in marking rubrics. (ii) No exercise asked the pupil to draw a machine. (iii) Two public chapters opened on the relation before the machine (\emph{cinqui\`eme}, order RMMM; MYP~5, order RMF, a relation and a set of pairs first), and four had no machine-tagged blocks (\emph{quatri\`eme}, \emph{seconde}, MYP~1, MYP~2), although the \emph{quatri\`eme} chapter served calculation programs run backwards. (iv) The \emph{troisi\`eme} chapter had no exercise on the uniqueness of the output, while five such written exercises existed at Grades~7 and~8: four were served only in MYP~5, and one in no chapter. (v) Written blocks were linked by no public chapter: the Grade~6 blocks on the machine's graph and their exercises, the opening problem and the definition of Grade~7. (vi) The robot was drawn for pupils of 14 to 17: in the opening problem of Grade~9, served in \emph{troisi\`eme} and MYP~5, and in the definition of Grade~11, served in DP~1.

\subsection{Revisions made after the audit}\label{sec:revisions}
The revisions were made on 2~November~2026 and form the state analysed as ``after''. In the source documents: the remark that draws the robot and the box side by side and names the box (546182); the definition of Grade~9 with the box and the line side by side and ``the inputs that the machine accepts'' (30614); boxes instead of robots in the opening problem of Grade~9 (30613); the robot removed from the definition of Grade~11 (4492); the machine and the sentence of covariation in ``in terms of'' (331804); the conversion drawn as a machine run both ways (444323); the calculation program as machines in a chain, run backwards in the reverse order (299656); ``the inputs that the machine accepts'' in the domain (3961); composition as two machines in a chain (3969) and the inverse introduced by $+3$ then $-3$ (3973); a written group of five exercises on drawing functions as machines at Grade~7 (546183--546191); the square machine run backwards at Grade~9 (546193); and repairs of existing blocks --- only the points are results of the machine (202177), joining the points of a table is a conjecture (30470, 422895), the situations of a class (291898), ``a rule that fits'' in the mystery machines (291818, 291826, 291844), and, in the two-rule machine, the question asking why the two machines give the same output for every input (291836). In the public chapters: the \emph{cinqui\`eme} chapter opens on the machine, the robot appearing once, and its eight exercises with the robot give way to the written exercises on the machine; the \emph{troisi\`eme} chapter receives the five written exercises on uniqueness, the drawing group and the square machine; the MYP~5 chapter presents the machine before the relation and the set of pairs. On 3~November two corrections followed: the coordinate plane of the \emph{troisi\`eme} chapter moved from the opening of the chapter to just before the graph, and the marked point of the vertical line test (Figure~\ref{fig:vertical}) was moved onto its curve. Tables~\ref{tab:corpus} and~\ref{tab:projections} give the counts before and after. What remains: the \emph{quatri\`eme}, MYP~1 and MYP~2 chapters have the machine in their lessons but in none of their exercises; the machine is still absent from the exercises of Grade~8 and Grade~10 beyond one block each; the Grade~6 graph definition (202177) and the Grade~7 opening problem and definition remain unlinked; exercise 291850 is now served in \emph{cinqui\`eme}.

\subsection{What the progression does not address}
Covariation is defined in the relation levels and said once at the joint, not practised in the function levels. Functions between sets other than numbers appear in a few exercises only. Functions given only by a graph are read, not constructed. The codomain is named at Grade~11; its role in the equality of two functions is made explicit in the first post-secondary year. The families of functions, which carry most of the later work on graphs and variation, are outside this paper.

\subsection{Hypotheses for an empirical evaluation}
\begin{hypothesis}
Pupils who met the constant machine and the two-rule machine at Grades 4--5 accept a constant function, and two formulas for one function, more often at 14--15 than pupils of a course that starts from the relation, on tasks of the kind used by \citet{Markovits1986}.
\end{hypothesis}
\begin{hypothesis}
Pupils who drew functions as machines in a chain compute $f(x+a)$ by substitution more often than pupils who did not, on the task reported by \citet{Oehrtman2008}.
\end{hypothesis}
\begin{hypothesis}
Compared with otherwise comparable pupils who did not meet this sentence, pupils who met the sentence on joining points from Grade~4 more often distinguish measured values from interpolated estimates, while still evaluating functions at intermediate inputs when their rules are known.
\end{hypothesis}
\begin{hypothesis}
On matched tasks, pupils who ran the square machine backwards more often correctly distinguish finding preimages from inverting a function than otherwise comparable pupils who did not.
\end{hypothesis}
No hypothesis is stated on covariation: the progression does not give it enough room in the function levels for a prediction.

\subsection{Limitations}
The author designed the progression, wrote the revisions and made the audit. Tags are regular expressions on the English halves; a tag records the presence of a marker, and every statement about a use was made by reading, but the reading was not double-coded. Tags include the marking rubrics embedded in exercise blocks, not only pupil-facing prompts, and they can miss relevant tasks: the drawing instruction of block 546185 is not tagged. The projections describe what the chapters show, not what is taught. The literature on the machine itself is thin and partly concerns adults \citep{McGowenDeMaroisTall2000}; the claims made from it are design hypotheses.

\section{Conclusion}\label{sec:conclusion}

One machine carries the progression from Grade~4 to Grade~10: an input, a rule, exactly one output. Its definition does not change; its drawing does, from the robot to the box to the arrow, each change shown once. On that machine the progression poses, from the age of nine, what makes the concept difficult: one input never has two outputs, two inputs may share one, a machine may answer the same number to everything, a rule need not be a formula, two rules may define one function, and a machine may be run backwards. The relation is kept for two quantities that vary together, and the set of pairs comes back as the graph. The analysis also locates where one machine is not enough: functions that are not numerical, functions given only by a graph, the domain and the codomain, and covariation. The audit showed that the machine had faded from the exercises once $f(x)$ was written, and the revisions brought it back in the lessons of the audited function and relation chapters and in the exercises of the \emph{cinqui\`eme} and \emph{troisi\`eme}; whether pupils taught this way understand functions differently is a question for an evaluation that has not been made.

\section*{Competing interests}
The author designed the progression analysed here and develops the platform on which it is published.

\section*{Acknowledgements}
The author used AI language tools while preparing this paper: Claude (Anthropic) assisted with the audit scripts, the drafting of the revisions and of the paper, and the search for and checking of bibliographic references, and ChatGPT (OpenAI) assisted with the critical reading of the plan, of the revisions and of the teacher's guide. The author designed the progression and takes full responsibility for the content.

\bibliographystyle{abbrvnat}
\bibliography{fonctions_refs}

\appendix
\section{Public resources}\label{app:resources}

The chapters cited are public at \url{https://commeunjeu.org}: French curriculum, \url{https://commeunjeu.org/cours/fr-grade-4-organizing-and-handling-data/}, \url{https://commeunjeu.org/cours/fr-grade-5-data-and-charts/}, \url{https://commeunjeu.org/cours/fr-grade-7-functions/}, \url{https://commeunjeu.org/cours/fr-grade-7-equations/}, \url{https://commeunjeu.org/cours/fr-grade-8-functions/}, \url{https://commeunjeu.org/cours/fr-grade-9-functions/}, \url{https://commeunjeu.org/cours/fr-grade-10-functions/}; International Baccalaureate, \url{https://commeunjeu.org/cours/ib-pyp-4-patterns-and-functions/}, \url{https://commeunjeu.org/cours/ib-pyp-4-data-handling/}, \url{https://commeunjeu.org/cours/ib-pyp-5-number-sequences-and-rules/}, \url{https://commeunjeu.org/cours/ib-myp-1-line-graphs/}, \url{https://commeunjeu.org/cours/ib-myp-2-coordinates-and-straight-line-graphs/}, \url{https://commeunjeu.org/cours/ib-myp-2-formulae/}, \url{https://commeunjeu.org/cours/ib-myp-5-functions/}, \url{https://commeunjeu.org/cours/ib-myp-5-extended-further-functions/}. A guide for teachers drawn from the same materials is at \url{https://commeunjeu.org/enseignants/enseigner-les-fonctions/}.

\section{Audit conventions}\label{app:audit}

The audit reads the English half of every block of the source documents (lessons and exercise sets of Grades 4 to 11 and of the two relation levels) at a given commit. A block is \emph{owned} by a file when it has an identifier and no link; a \emph{mirror} carries a link and a copy of the body, and is counted apart. A block is tagged \emph{machine} when its text contains the word ``machine'' or the image of the machine; \emph{draw} when its text has ``draw'', ``complete'' or ``fill in'' followed within 40 characters by ``machine'', which misses some drawing instructions; \emph{not a function} when it presents a correspondence as not a function, or names the vertical line test, two outputs, or more than one output; \emph{constant} when it names the constant function; \emph{formula} when it contains an equality of the form $f(x) = \ldots$. Answers are counted apart from exercises. The projections script follows each link of a curriculum chapter to its source block and applies the machine tag to that block; it records the order of the lesson's blocks that contain the machine (M), define a relation (R), or define a function without the machine (F). The two states are the commits \texttt{596ec5ecec} (before) and \texttt{e8927be52e} (after) for the source documents, and \texttt{596ec5ecec} and \texttt{d809e0071a} for the projections.

\end{document}